\documentclass[12pt]{article}

\usepackage{amsmath,amssymb,amsbsy,amsfonts,amsthm,latexsym,
                        amsopn,amstext,amsxtra,euscript,amscd}

\newtheorem{theorem}{Theorem}[section]
\newtheorem{lemma}[theorem]{Lemma}

\newtheorem{corollary}[theorem]{Corollary}

\theoremstyle{definition}

\def\fl#1{\left\lfloor#1\right\rfloor}
\def\rf#1{\left\lceil#1\right\rceil}

\def\N{\mathbb{N}}

\def\cB{\mathcal B}

\def\cI{\mathcal I}

\def\cM{\mathcal M}
\def\cN{\mathcal N}

\def\Z{{\mathbb Z}}
\def\R{{\mathbb R}}

\def\eps{{\varepsilon}}
\def\e{{\rm\bf e\/}}

\def\mand{\qquad \mbox{and} \qquad}

\def\\{\cr}
\def\({\left(}
\def\){\right)}
\def\[{\left[}
\def\]{\right]}
\def\<{\langle}
\def\>{\rangle}
\def\fl#1{\left\lfloor#1\right\rfloor}
\def\rf#1{\left\lceil#1\right\rceil}
\def\le{\leqslant}
\def\ge{\geqslant}

\begin{document}

\title{\sc Prime numbers with\break Beatty sequences}

\author{
{\sc William D.~Banks} \\
{Department of Mathematics} \\
{University of Missouri} \\
{Columbia, MO 65211 USA} \\
{\tt bbanks@math.missouri.edu} \\
\and
{\sc Igor E.~Shparlinski} \\
{Department of Computing}\\
{Macquarie University} \\
{Sydney, NSW 2109, Australia} \\
{\tt igor@ics.mq.edu.au}}

\maketitle
\begin{abstract}
A study of certain Hamiltonian systems has lead Y.~Long to
conjecture the existence of infinitely many primes of the form
$p=2\fl{\alpha n}+1$, where $1<\alpha<2$ is a fixed irrational
number. An argument of P.~Ribenboim coupled with classical results
about the distribution of fractional parts of irrational multiples
of primes in an arithmetic progression immediately imply that this
conjecture holds in a much more precise asymptotic form. Motivated
by this observation, we give an asymptotic formula for the number
of primes $p=q\fl{\alpha n+\beta}+a$ with $n\le N$, where
$\alpha,\beta$ are real numbers such that $\alpha$ is positive and
irrational of finite type (which is true for almost all $\alpha$)
and $a,q$ are integers with $0\le a<q\le N^\kappa$ and
$\gcd(a,q)=1$, where $\kappa>0$ depends only on $\alpha$. We also
prove a similar result for primes $p=\fl{\alpha n+\beta}$ such
that $p\equiv a\pmod q$.
\end{abstract}

\newpage

\section{Introduction}

For two fixed real numbers $\alpha$ and $\beta$, the corresponding
\emph{non-homogeneous Beatty sequence} is the sequence of integers
defined by
$$
\cB_{\alpha,\beta}=\(\fl{\alpha n+\beta}\)_{n=1}^\infty.
$$
Beatty sequences appear in a variety of apparently unrelated
mathematical settings, and because of their versatility, the
arithmetic properties of these sequences have been extensively
explored in the literature; see, for example,
\cite{Abe,BaSh1,BaSh2,Beg,FraHolz,Kom1,Kom2,LuZh,OB,Tijd} and the
references contained therein.

In 2000, while investigating the Maslov-type index theory for
Hamiltonian systems, Long~\cite{Long} made the following
conjecture:

\begin{quote}
{\sc Conjecture.} \emph{For every irrational number $1<\alpha<2$,
there are infinitely many prime numbers of the form $p=2\fl{\alpha
n}+1$ for some $n\in\N$.}
\end{quote}

\noindent Jia~\cite{Jia} has given a lower bound for the number of
such primes $p$ in the interval $(x/2,x]$.  We remark that, using
a simple modification to an argument given by
Ribenboim~\cite[Chapter~4.V]{Riben}, one can show further that the
number of such primes $p\le x$ is asymptotic to
$\alpha^{-1}\pi(x)$ as $x\to\infty$; see also~\cite{LW}. Moreover,
Ribenboim's method also applies to the general problem of
estimating
$$
\cN_{\alpha,\beta;q,a}(x)=\#\big\{n\le x~:~p=q\fl{\alpha
n+\beta}+a\text{~is prime}\big\},
$$
where $\alpha,\beta$ are fixed real numbers such that $\alpha$ is
positive and irrational, and $a,q$ are integers with $0\le a<q$
and $\gcd(a,q)=1$. In fact, if $a$ and $q$ are fixed, one easily
derives the asymptotic formula
$$
\cN_{\alpha,\beta;q,a}(x)=(1+o(1))\,
\frac{q}{\varphi(q)}\,\pi(x)\qquad(x\to\infty),
$$
where the function implied by $o(\cdot)$ depends on $\alpha$,
$\beta$ and $q$, and $\varphi(\cdot)$ is the \emph{Euler
function}. Motivated by this observation, we consider here the
problem of finding uniform estimates for
$\cN_{\alpha,\beta;q,a}(x)$ if $q$ is allowed to grow with $x$. We
also consider the same problem for the counting function
$$
\cM_{\alpha,\beta;q,a}(x)=\#\big\{n\le x~:~p=\fl{\alpha
n+\beta}\text{~is prime, and~}p\equiv a\pmod q\big\}.
$$
In particular, in the case that $\alpha$ is of \emph{finite type}
(which is true for \emph{almost all} $\alpha$ in sense of Lebesgue
measure), our main results yield (by partial summation) nontrivial
results for both $\cN_{\alpha,\beta;q,a}(x)$ and
$\cM_{\alpha,\beta;q,a}(x)$ even if $q$ grows as a certain power
of $x$.

\bigskip

\noindent\textbf{Acknowledgements.} This work started in August of
2006, while one of the authors (I.~S.) attended  the {\it Fourth
China-Japan Conference on Number Theory} at Shandong University in
Weihai, China; I.~S.\ thanks the organizers, Shigeru Kanemitsu and
Jianya Liu, for their hospitality
and the opportunity of participating in this event. The authors
would like to thank Wesley Nevans for suggesting a generalization
of the question we had originally intended to study, Christian
Mauduit for calling our attention to the result of
Ribenboim~\cite{Riben} on primes in a Beatty sequence, and Ahmet
G\"ulo\u glu for pointing out a mistake in the original version of
the manuscript. During the preparation of this paper, I.~S.\ was
supported in part by ARC grant DP0556431.

\section{Notation}

The notation $\|x\|$ is used to denote the distance from the real
number $x$ to the nearest integer; that is,
$$
\|x\|=\min_{n\in\Z}|x-n|\qquad(x\in\R).
$$
As usual, we denote by $\fl{x}$, $\rf{x}$, and $\{x\}$ the
greatest integer $\le x$, the least integer $\ge x$, and the
fractional part of $x$, respectively.

We also put $\e(x)=e^{2\pi ix}$ for all real numbers $x$ and use
$\Lambda(\cdot)$ to denote the \emph{von Mangoldt function}:
$$
\Lambda(n)=
\begin{cases}
\log p&\quad\text{if $n$ is a power of a prime $p$;} \\
0&\quad\text{otherwise.}
\end{cases}
$$

Throughout the paper, the implied constants in symbols $O$, $\ll$
and  $\gg$ may depend on the parameters $\alpha$ and $\beta$ but
are absolute unless indicated otherwise. We recall that the
notations $A\ll B$, $B\gg A$ and $A=O(B)$ are all equivalent to
the statement that $|A|\le c|B|$ for some constant $c>0$.

\section{Preliminaries}
\label{sec:prelims}

Recall that the \emph{discrepancy} $D(M)$ of a sequence of (not
necessarily distinct) real numbers $a_1,a_2,\ldots,a_M\in[0,1)$ is
defined by
\begin{equation}
\label{eq:descr_defn}
D(M)=\sup_{\cI\subseteq[0,1)}\left|\frac{V(\cI,M)}{M}-|\cI|\,\right|,
\end{equation}
where the supremum is taken all subintervals $\cI=(c,d)$ of the
interval $[0, 1)$, $V(\cI,M)$ is the number of positive integers
$m\le M$ such that $a_m\in\cI$, and $|\cI|=d-c$ is the length of
$\cI$.

For an irrational number $\gamma$, we define its \emph{type}
$\tau$ by the relation
$$
\tau=\sup\Bigl\{\varrho \in\R~:~\liminf_{\substack{n\to\infty\\
(n\in\N)}} ~n^\varrho\,\|\gamma n\|=0\Bigr\}.
$$
Using~\emph{Dirichlet's approximation theorem}, it is easily seen
that $\tau\ge 1$ for every irrational number $\gamma$. The
celebrated theorems of Khinchin~\cite{Khin} and of
Roth~\cite{Roth1, Roth2} assert that $\tau=1$ for almost all real
(in the sense of the Lebesgue measure) and all irrational
algebraic numbers $\gamma$, respectively; see
also~\cite{Bug,Schm}.

For every irrational number $\gamma$, it is well known that the
sequence of fractional parts
$\{\gamma\},\{2\gamma\},\{3\gamma\},\,\ldots\,,$ is
\emph{uniformly distributed modulo~$1$} (for instance, see
\cite[Example~2.1, Chapter~1]{KuNi}). If $\gamma$ is of finite
type, this statement can be made more precise.
Let $D_{\gamma,\delta}(M)$ denote the discrepancy of the sequence
of fractional parts $(\{\gamma m + \delta\})_{m=1}^M$.
By~\cite[Theorem~3.2, Chapter~2]{KuNi} we have:

\begin{lemma}
\label{lem:discr with type}  Let $\gamma$ be a fixed irrational
number of finite type $\tau<\infty$.  Then, for all $\delta\in\R$
the following bound holds:
$$
D_{\gamma,\delta}(M)\le M^{-1/\tau+o(1)}\qquad(M\to\infty),
$$
where the function implied by $o(\cdot)$ depends only on $\gamma$.
\end{lemma}

The following elementary result characterizes the set of numbers
that occur in a Beatty sequence $\cB_{\alpha,\beta}$ in the case
that $\alpha>1$:

\begin{lemma}
\label{lem:Beatty values} Let $\alpha,\beta\in\R$ with $\alpha >
1$. Then, an integer $m$ has the form $m= \fl{\alpha n + \beta}$
for some integer $n$ if and only if
$$
0<\bigl\{\alpha^{-1}(m-\beta+1)\bigr\}\le\alpha^{-1}.
$$
The value of $n$ is determined uniquely by $m$.
\end{lemma}

\begin{proof}
It is easy to see that an integer $m$ has the form $m=\fl{\alpha
n+\beta}$ for some integer $n$ if and only if the inequalities
$$
\frac{m-\beta}{\alpha} \le n < \frac{m-\beta+1}{\alpha}
$$
hold, and since $\alpha>1$ the value of $n$ is determined
uniquely.
\end{proof}

We also need the following statement, which is a simplified and
weakened version of a theorem of Balog and Perelli~\cite{BalPer}
(see also~\cite{Lav}):

\begin{lemma}
\label{lem:balog_perelli} For an arbitrary real number $\vartheta
$ and coprime integers $a,q$ with $0\le a<q$, if $|\vartheta
-b/d|\le 1/L$ and $\gcd(b,d)=1$, then the bound
$$
\sum_{\substack{n\le L\\n\equiv a\pmod q}}\Lambda(n)\,\e(\vartheta
n) \ll\(\frac{L}{d^{1/2}}+d^{1/2}L^{1/2}+L^{4/5}\)(\log L)^3
$$
holds, where the implied constant is absolute.
\end{lemma}

Finally, we use the \emph{Siegel--Walfisz theorem} (see, for
example, the book~\cite{Hux} by Huxley), which asserts that for
any fixed constant $B>0$ and uniformly for integers $L\ge 3$ and
$0\le a<q\le(\log L)^B$ with $\gcd(a,q)=1$, one has
\begin{equation}
\label{eq:siegelwalfisz} \sum_{\substack{n\le L\\n\equiv a\pmod
q}}\Lambda(n)=\frac{L}{\varphi(q)}+O\(L\,\exp\(-C_B\sqrt{\log
L}\,\)\),
\end{equation}
where $C_B>0$ is an absolute constant that depends only on $B$.

\section{Bounds on exponential sums}

The following result may be well known but does not seem to be
recorded in the literature. Thus, we present it here with a
complete proof.

\begin{theorem}
\label{thm:non-liou-theta} Let $\gamma$ be a fixed irrational
number of finite type $\tau<\infty$.  Then, for every real number
$0<\eps<1/(8\tau)$, there is a number $\eta>0$ such that the bound
$$
\left|\sum_{m\le M}\Lambda(qm+a)\,\e(\gamma k m)\right|\le
M^{1-\eta}
$$
holds for all integers $1\le k\le M^\eps$ and $0\le a<q\le
M^{\eps/4}$ with $\gcd(a,q)=1$ provided that $M$ is sufficiently
large.
\end{theorem}

\begin{proof}
Fix a constant $\varrho$ such that
\begin{equation}
\label{eq:eps cond} 1\le \tau<\varrho<\frac{1}{8\eps}
\end{equation}
Since $\gamma$ is of type $\tau$, for some constant $c>0$ we have
\begin{equation}
\label{eq:Liou} \|\gamma d\|>c d^{-\varrho}\qquad(d\ge 1).
\end{equation}

Let $k,a,q$ be integers with the properties stated in the
proposition, and write
\begin{equation}
\label{eq:lamblamb} \sum_{m\le M}\Lambda(qm+a)\,\e(\gamma k
m)=\e(-\vartheta a)\sum_{\substack{n\le L\\n\equiv a\pmod
q}}\Lambda(n)\,\e(\vartheta n),
\end{equation}
where $\vartheta =\gamma k/q$ and $L=qM+a$. Let $b/d$ be the
convergent in the continued fraction expansion of $\vartheta $
which has the largest denominator $d$ not exceeding $L^{1-\eps}$;
then,
\begin{equation}
\label{eq:theta-mix1} \left|\frac{\gamma k}{q}-\frac{b}{d}\right|
\le\frac{1}{d L^{1-\eps}}.
\end{equation}
Multiplying by $qd$, we get from~\eqref{eq:Liou}:
$$
\frac{q}{L^{1-\eps}}\ge\left|\gamma kd-bq\right| \ge\|\gamma
kd\|>c (kd)^{-\varrho}.
$$
Thus, since $k\le L^\eps$ and $q\le L^{\eps/4} \le L^\eps$,  we
see that under the condition~\eqref{eq:eps cond} the bound
\begin{equation}
\label{eq:theta-mix2} d \ge  C L^{(1-2 \eps)/\varrho - \eps} \ge C
L^{1/(4\varrho)}
\end{equation}
holds, where $C=c^{1/\varrho}$ and $L$ is sufficiently large.

Inserting~\eqref{eq:theta-mix2} into~\eqref{eq:theta-mix1} and
using~\eqref{eq:eps cond} again, we conclude  that
$$
\left|\frac{\gamma k}{q}-\frac{b}{d}\right| \le\frac{1}{ C L^{1+
1/(4\varrho)-\eps}}\le\frac{1}{L}
$$
if $L$ is sufficiently large. We are therefore in a position to
apply Lemma~\ref{lem:balog_perelli}; taking into
account~\eqref{eq:eps cond}, \eqref{eq:theta-mix2}, and the fact
that $d\le L^{1-\eps}$, it follows that the bound
$$
\sum_{\substack{n\le L\\n\equiv a\pmod q}}\Lambda(n)\,\e(\vartheta
n) \ll \(L^{1-1/(8\varrho)} + L^{1-\eps/2}\)(\log L)^3 \le
L^{1-\eps/3}
$$
holds for all sufficiently large $L$. Since $L\ll qM\le
M^{1+\eps/4}$, the result now follows from simple calculations
after inserting this estimate into~\eqref{eq:lamblamb}.
\end{proof}

Using similar arguments, we have:

\begin{theorem}
\label{thm:non-liou-theta2} Let $\gamma$ be a fixed irrational
number of finite type $\tau<\infty$.  Then, for every real number
$0<\eps<1/(8\tau)$, there is a number $\eta>0$ such that the bound
$$
\left|\sum_{\substack{m\le M\\m\equiv a\pmod
q}}\Lambda(m)\,\e(\gamma k m)\right|\le M^{1-\eta}
$$
holds for all integers $1\le k\le M^\eps$ and $0\le a<q\le
M^{\eps/4}$ with $\gcd(a,q)=1$ provided that $M$ is sufficiently
large.
\end{theorem}

\section{Main Results}

\begin{theorem}
\label{thm:main} Let $\alpha$ and $\beta$ be a fixed real numbers
with $\alpha$ positive, irrational, and of finite type. Then there
is a positive constant $\kappa>0$ such that for all integers $0\le
a<q\le N^\kappa $ with $\gcd(a,q)=1$, we have
$$
\sum_{n\le N}\Lambda(q\fl{\alpha n+\beta}+a)=\alpha^{-1}\sum_{m\le
\fl{\alpha N+\beta}}\Lambda(qm+a)+O\(N^{1-\kappa}\)
$$
where the implied constant depends only on $\alpha$ and $\beta$.
\end{theorem}

\begin{proof}
Suppose first that $\alpha>1$. It is obvious that if $\alpha$ is
of finite type, then so is $\alpha^{-1}$. We choose
$$
0<\eps<\frac{1}{16\tau},
$$
where  $1\le\tau<\infty$ is the type of $\alpha^{-1}$.

First, let us suppose that $\alpha>1$. Put $\gamma=\alpha^{-1}$,
$\delta = \alpha^{-1}(1-\beta)$, and $M=\fl{\alpha N+\beta}$. By
Lemma~\ref{lem:Beatty values}, it follows that
\begin{equation}
\label{eq:S}
\begin{split}
S_{\alpha,\beta;q,a}(N)&=\sum_{n\le N}\Lambda(q\fl{\alpha
n+\beta}+a)\\
&=\sum_{\substack{m\le M\\0<\{\gamma m+\delta\}\le\gamma}}
\Lambda(qm+a)+O(1)\\
&=\sum_{m\le M} \Lambda(qm+a)\,\psi(\gamma m+\delta)+O(1),
\end{split}
\end{equation}
where $\psi(x)$ is the periodic function with period one for which
$$
\psi(x) = \left\{  \begin{array}{ll}
1& \quad \hbox{if $0< x\le \gamma$}; \\
0& \quad \mbox{if $\gamma<x\le 1$}.
\end{array} \right.
$$
By a classical result of Vinogradov (see~\cite[Chapter~I,
Lemma~12]{Vin}) it is known that for any $\Delta$ such that
$$
0 < \Delta < \frac{1}{8} \mand \Delta\le
\frac{1}{2}\min\{\gamma,1-\gamma\},
$$
there is a real-valued function $\psi_\Delta(x)$ with the
following properties:
\begin{itemize}
\item[$(i)$~~] $\psi_\Delta(x)$ is periodic with period one;

\item[$(ii)$~~] $0 \le \psi_\Delta(x) \le 1$ for all $x\in\R$;

\item[$(iii)$~~] $\psi_\Delta(x) = \psi(x)$ if $\Delta\le x\le
\gamma-\Delta$ or if $\gamma+\Delta\le x\le 1-\Delta$;

\item[$(iv)$~~] $\psi_\Delta(x)$ can be represented as a Fourier
series
$$
\psi_\Delta(x) = \gamma + \sum_{k=1}^\infty \bigl(\,g_k\,\e(kx) +
h_k \,\e(-kx)\bigr),
$$
where the coefficients satisfy the uniform bound
\begin{equation}
\label{eq:coeffbounds} \max\bigl\{|g_k|, |h_k|\bigr\}
\ll\min\bigl\{k^{-1},k^{-2}\Delta^{-1}\bigr\} \qquad(k\ge 1).
\end{equation}
\end{itemize}
Therefore, from~\eqref{eq:S} we deduce that
\begin{equation}
\label{eq:Char Fun Approx} S_{\alpha,\beta;q,a}(N) = \sum_{m\le
M}\Lambda(qm+a)\,\psi_\Delta(\gamma m + \delta)
+O\bigl(1+V(\cI,M)\log N\bigr),
\end{equation}
where $V(\cI,M)$ denotes the number of positive integers $m\le M$
such that
$$
\{\gamma m+\delta\}\in\cI=
[0,\Delta)\cup(\gamma-\Delta,\gamma+\Delta) \cup(1-\Delta,1).
$$
Since $|\cI|\ll\Delta$, it follows from the
definition~\eqref{eq:descr_defn} and Lemma~\ref{lem:discr with
type}, that
\begin{equation}
\label{eq:bound V(I,M) with type} V(\cI,M)\ll\Delta N+N^{1-\eps},
\end{equation}
where the implied constant depends only on $\alpha$.

To estimate the sum in~\eqref{eq:Char Fun Approx}, we use the
Fourier expansion for $\psi_\Delta(\gamma m + \delta)$ and change
the order of summation, obtaining
\begin{equation}
\label{eq:triplesplit}
\begin{split}
\sum_{m\le M}\Lambda(qm+a)\,\psi_\Delta(\gamma m &+ \delta)\\
=\gamma\sum_{m\le M}\Lambda(qm+a) &+\sum_{k=1}^\infty
g_k\,\e(\delta k)
\sum_{m\le M}\Lambda(qm+a)\,\e(\gamma k m) \\
&+\sum_{k=1}^\infty h_k\,\e(-\delta k) \sum_{m\le
M}\Lambda(qm+a)\,\e(-\gamma k m).
\end{split}
\end{equation}

By Theorem~\ref{thm:non-liou-theta} and the
bound~\eqref{eq:coeffbounds}, we see that for $0\le a<q\le
M^{\eps/4}$, we have
\begin{equation}
\label{eq:bd5} \sum_{k\le M^\eps} g_k\,\e(\delta k) \sum_{m\le
M}\Lambda(qm+a)\,\e(\gamma k m) \ll M^{1-\eta} \sum_{k\le
M^{\eps}} k^{-1} \ll M^{1-\eta/2} ,
\end{equation}
for some $\eta>0$ that depends only on $\alpha$. Similarly,
\begin{equation}
\label{eq:bd6} \sum_{k\le M^{\eps}} h_k\,\e(-\delta k) \sum_{m\le
M}\Lambda(qm+a)\,\e(-\gamma k m) \ll M^{1-\eta/2}.
\end{equation}
On the other hand, using the trivial bound
$$
\left|\sum_{m\le M}\Lambda(qm+a)\,\e(\gamma k m)\right|
\le\sum_{n\le N}\Lambda(n)\ll N,
$$
we have
\begin{equation}
\label{eq:bd7} \sum_{k>M^{\eps}} g_k\,\e(\delta k) \sum_{m\le
M}\Lambda(qm+a)\,\e(\gamma k m) \ll N\sum_{k>M^{\eps}}
k^{-2}\Delta^{-1}\ll N^{1-\eps}\Delta^{-1},
\end{equation}
and
\begin{equation}
\label{eq:bd8} \sum_{k>M^{\eps}} h_k\,\e(-\delta k) \sum_{m\le
M}\Lambda(qm+a)\,\e(-\gamma k m) \ll N^{1-\eps}\Delta^{-1}.
\end{equation}
Inserting the bounds, \eqref{eq:bd5}, \eqref{eq:bd6},
\eqref{eq:bd7} and \eqref{eq:bd8} into~\eqref{eq:triplesplit}, we
obtain
\begin{equation}
\begin{split}
\label{eq:main term}
    \sum_{m\le M}\Lambda(qm+a)\,&\psi_\Delta(\gamma m + \delta) \\
=&\gamma\sum_{m\le M}\Lambda(qm+a)  + O(M^{1-\eta/2} +
N^{1-\eps}\Delta^{-1}),
\end{split}
\end{equation}
where the constant implied by $O(\cdot)$ depends only on $\alpha$
and $\beta$.

Substituting~\eqref{eq:bound V(I,M) with type} and~\eqref{eq:main
term} in~\eqref{eq:Char Fun Approx} and choosing
$\Delta=N^{-\eps/4}$, it follows that
\begin{equation}
\label{eq:ryan} S_{\alpha,\beta;q,a}(N)=\gamma \sum_{m\le
M}\Lambda(qm+a)+O\(N^{1-\kappa}\),
\end{equation}
for some $\kappa$ which depends only on $\alpha$. This concludes
the proof  in the case that $\alpha>1$.

If $\alpha<1$, we put $t=\rf{\alpha^{-1}}$ and write
$$
\sum_{n\le N}\Lambda(q\fl{\alpha n+\beta}+a)= \sum_{j=0}^{t-1}
\sum_{n\le (N-j)/t} \Lambda(q\fl{\alpha tn+\alpha j+\beta}+a).
$$
Applying the preceding argument with the irrational number $\alpha
t>1$, we conclude the proof.
\end{proof}

In particular, using the Siegel--Walfisz
theorem~\eqref{eq:siegelwalfisz} to estimate the sum
in~\eqref{eq:ryan} for ``small'' $a$ and $q$, we obtain:

\begin{corollary}
\label{cor:good-theta 1} Under the conditions of
Theorem~\ref{thm:main}, for any constant $B>0$ and uniformly for
all integers $N\ge 3$ and $0\le a<q\le(\log N)^B$ with
$\gcd(a,q)=1$, we have
$$
\sum_{n\le N}\Lambda(q\fl{\alpha
n+\beta}+a)=\frac{q}{\varphi(q)}\,N+O\(N\,\exp\(-C\sqrt{\log
N}\,\)\)
$$
for some constant $C>0$ that depends only on $\alpha$, $\beta$ and
$B$.
\end{corollary}

In the special case that $(a,q) = (0,1)$ or $(1,2)$ (the latter
case corresponding to primes in the Long conjecture), we can use a
well known bound on the error term in the Prime Number Theorem
(proved independently by Korobov~\cite{Kor} and
Vinogradov~\cite{Vin2}) to achieve the following sharper result:

\begin{corollary}
\label{cor:good-theta 2} Suppose that $(a,q)= (0,1)$ or
$(a,q)=(1,2)$. Then, under the conditions of
Theorem~\ref{thm:main}, for any constant $B>0$ and uniformly for
all integers $N\ge 3$, we have
$$
\sum_{n\le N}\Lambda(q\fl{\alpha
n+\beta}+a)=q\,N+O\(N\,\exp\(-c(\log N)^{3/5}(\log\log
N)^{-1/5}\)\)
$$
for some absolute constant $c>0$.
\end{corollary}

Finally, using Lemma~\ref{thm:non-liou-theta2} in place of
Lemma~\ref{thm:non-liou-theta}, we obtain the following analogues
of Theorem~\ref{thm:main} and its two corollaries:

\begin{theorem}
\label{thm:main2} Let $\alpha$ and $\beta$ be a fixed real numbers
with $\alpha$ positive, irrational, and of finite type. Then there
is a positive constant $\kappa>0$ such that for all integers $0\le
a<q\le N^\kappa $ with $\gcd(a,q)=1$, we have
$$
\sum_{\substack{n\le N\\\fl{\alpha n+\beta}\equiv a\pmod
q}}\Lambda(\fl{\alpha n+\beta})=\alpha^{-1}\sum_{\substack{m\le
\fl{\alpha N+\beta}\\m\equiv a\pmod
q}}\Lambda(m)+O\(N^{1-\kappa}\)
$$
where the implied constant depends only on $\alpha$ and $\beta$.
\end{theorem}

\begin{corollary}
\label{cor:good-theta 3} Under the conditions of
Theorem~\ref{thm:main2}, for any constant $B>0$ and uniformly for
all integers $N\ge 3$ and $0\le a<q\le(\log N)^B$ with
$\gcd(a,q)=1$, we have
$$
\sum_{\substack{n\le N\\\fl{\alpha n+\beta}\equiv a\pmod
q}}\Lambda(\fl{\alpha
n+\beta})=\frac{N}{\varphi(q)}+O\(N\,\exp\(-C\sqrt{\log N}\,\)\)
$$
for some constant $C>0$ that depends only on $\alpha$, $\beta$ and
$B$.
\end{corollary}

\begin{corollary}
\label{cor:good-theta 4} Suppose that $(a,q)= (0,1)$ or
$(a,q)=(1,2)$. Then, under the conditions of
Theorem~\ref{thm:main2}, for any constant $B>0$ and uniformly for
all integers $N\ge 3$, we have
$$
\sum_{\substack{n\le N\\\fl{\alpha n+\beta}\equiv a\pmod
q}}\Lambda(\fl{\alpha n+\beta})=N+O\(N\,\exp\(-c(\log
N)^{3/5}(\log\log N)^{-1/5}\)\)
$$
for some absolute constant $c>0$.
\end{corollary}

\end{document}